\documentclass[12pt]{amsart}
\usepackage{amsmath,amssymb,amsbsy,amsfonts,amsthm,latexsym,mathabx,
            amsopn,amstext,amsxtra,euscript,amscd,stmaryrd,mathrsfs,
            cite,array,mathtools,enumerate}
            
\usepackage{url}
\usepackage[colorlinks,linkcolor=blue,anchorcolor=blue,citecolor=blue,backref=page]{hyperref}

\usepackage{color}

\renewcommand*{\backref}[1]{}
\renewcommand*{\backrefalt}[4]{%
    \ifcase #1 (Not cited.)%
    \or        (p.\,#2)%
    \else      (pp.\,#2)%
    \fi}
\begin{document}

\newtheorem{theorem}{Theorem}
\newtheorem{lemma}[theorem]{Lemma}
\newtheorem{claim}[theorem]{Claim}
\newtheorem{cor}[theorem]{Corollary}
\newtheorem{prop}[theorem]{Proposition}
\newtheorem{definition}{Definition}
\newtheorem{question}[theorem]{Open Question}
\newtheorem{conj}[theorem]{Conjecture}
\newtheorem{prob}{Problem}
\newtheorem{algorithm}[theorem]{Algorithm}

\numberwithin{equation}{section}
\numberwithin{theorem}{section}

\def\squareforqed{\hbox{\rlap{$\sqcap$}$\sqcup$}}
\def\qed{\ifmmode\squareforqed\else{\unskip\nobreak\hfil
\penalty50\hskip1em\nobreak\hfil\squareforqed
\parfillskip=0pt\finalhyphendemerits=0\endgraf}\fi}

\def\cA{{\mathcal A}}
\def\cB{{\mathcal B}}
\def\cC{{\mathcal C}}
\def\cD{{\mathcal D}}
\def\cE{{\mathcal E}}
\def\cF{{\mathcal F}}
\def\cG{{\mathcal G}}
\def\cH{{\mathcal H}}
\def\cI{{\mathcal I}}
\def\cJ{{\mathcal J}}
\def\cK{{\mathcal K}}
\def\cL{{\mathcal L}}
\def\cM{{\mathcal M}}
\def\cN{{\mathcal N}}
\def\cO{{\mathcal O}}
\def\cP{{\mathcal P}}
\def\cQ{{\mathcal Q}}
\def\cR{{\mathcal R}}
\def\cS{{\mathcal S}}
\def\cT{{\mathcal T}}
\def\cU{{\mathcal U}}
\def\cV{{\mathcal V}}
\def\cW{{\mathcal W}}
\def\cX{{\mathcal X}}
\def\cY{{\mathcal Y}}
\def\cZ{{\mathcal Z}}

\def\fA{{\mathfrak A}}
\def\fJ{{\mathfrak J}}

\def\fE{{\mathfrak E}}
\def\fS{{\mathfrak S}}

\def\sssum{\mathop{\sum\!\sum\!\sum}}
\def\ssum{\mathop{\sum\ldots \sum}}

\def\Xm{\cX_m}

\def \C {{\mathbb C}}
\def \F {{\mathbb F}}
\def \L {{\mathbb L}}
\def \K {{\mathbb K}}
\def \N {{\mathbb N}}
\def \R {{\mathbb R}}
\def \Q {{\mathbb Q}}
\def \Z {{\mathbb Z}}

\def\sA{\mathscr A}
\def\sB{\mathscr B}
\def\sL{\mathscr L}

\def\barG{\overline{\cG}}
\def\\{\cr}
\def\({\left(}
\def\){\right)}
\def\fl#1{\left\lfloor#1\right\rfloor}
\def\rf#1{\left\lceil#1\right\rceil}

\newcommand{\pfrac}[2]{{\left(\frac{#1}{#2}\right)}}

\def\rem{\mathrm{\, rem~}}

\def \Prob{{\mathrm {}}}
\def\e{\mathbf{e}}
\def\ep{{\mathbf{\,e}}_p}
\def\eq{{\mathbf{\,e}}_{q}}
\def\er{{\mathbf{\,e}}_r}
\def\es{{\mathbf{\,e}}_s}
\def\eM{{\mathbf{\,e}}_M}
\def\eps{\varepsilon}
\def\Res{\mathrm{Res}}
\def\vec#1{\mathbf{#1}}

\def \li {\mathrm {li}\,}

\def\ip{\overline p}
\def\ipd{\ip_d}
\def\iq{\overline q}

\def\e{{\mathbf{\,e}}}
\def\ep{{\mathbf{\,e}}_p}
\def\em{{\mathbf{\,e}}_m}

\def\mand{\qquad\mbox{and}\qquad}

\newcommand{\commF}[1]{\marginpar{%
\begin{color}{red}
\vskip-\baselineskip 
\raggedright\footnotesize
\itshape\hrule \smallskip F: #1\par\smallskip\hrule\end{color}}}

\newcommand{\commI}[1]{\marginpar{%
\begin{color}{blue}
\vskip-\baselineskip 
\raggedright\footnotesize
\itshape\hrule \smallskip I: #1\par\smallskip\hrule\end{color}}}

\title[Character sums with smooth numbers]{Character sums with smooth numbers}
\date{\today}

\author{Igor E.  Shparlinski}

\address{Department of Pure Mathematics, University of New South Wales\\
2052 NSW, Australia.}

\email{igor.shparlinski@unsw.edu.au}


\begin{abstract} We use the large sieve inequality for smooth numbers due to 
S.~Drappeau, A.~Granville and X.~Shao (2017),  together with some other arguments, 
to improve  their bounds on the frequency of  pairs $(q,\chi)$  of moduli $q$ and primitive 
characters $\chi$ modulo $q$, for which the corresponding character sums with smooth numbers are large. 
\end{abstract}
 
\maketitle

\section{Introduction}

Let $\Psi(x,y)$ be the set of $y$-smooth integers $n \le x$, 
that is, 
$$
\Psi(x,y) = \{n \in \Z\cap[1,x]~:~ P(n) \le y\}, 
$$
where $P(n)$ is the largest prime divisor of a positive 
integer $n$. 
We also denote the cardinality of $\Psi(x,y)$ by   $\psi(x,y) = \# \Psi(x,y)$.  

Let $\cX_q$ denote the set of all $\varphi(q)$  multiplicative characters modulo an integer  $q\ge 2$
and let $\cX_q^*$ be the set of primitive characters $\chi\in \cX_q$,
where $\varphi(q)$ denotes the Euler function of $q$, 
we refer to~\cite[Chapter~3]{IwKow} for a background on characters.

Given a sequence  $\sA = \{a_n\}_{n=1}^\infty$ of complex numbers, 
w now consider the character sums
$$
S_q(\chi;\sA,x,y) =  \sum_{n \in \Psi(x,y)} a_n \chi(n).
$$
In the case when $a_n =1$, $n =1, 2, \ldots$, we simply write 
$$
S_q(\chi;x,y) =  \sum_{n \in \Psi(x,y)}  \chi(n).
$$
 Drappeau, Granville and Shao~\cite{DGS} have recently shown that there
  exist absolute constants $C_0, c_0>0$ such that  real 
$x$, $y$ and $Q$ satisfy
$$
(\log x)^{C_0} \le y \le x \mand Q \le \min \{y^{c_0}, \exp(c_0 \log x/\log \log x)\}, 
$$
then for any fixed $\beta\ge 0$ for all but at most 
\begin{equation}
\label{eq:bound E}
E = O\((\log x)^{3\beta+13}\)
\end{equation}
pairs $(q, \chi)$ with a positive integer $q \le Q$ and $\chi \in \cX_q^*$, 
for every $t \in [x^{1/4},x]$ we have 
\begin{equation}
\label{eq:small S}
\left|S_q(\chi;t,y)\right|  < \frac{\psi(t,y)}{(u \log u)^4 (\log x)^\beta},
\end{equation}
where as usual
$$
u =  \frac{\log x}{\log y}. 
$$
The result is based on the obtained in~\cite{DGS} large sieve inequality 
with smooth numbers. It seems that the same result also holds for the 
sums $S_q(\chi;\sA, t,y)$ satisfying the bound~\eqref{eq:small S}. 

Here we modify the scheme of the proof from~\cite{DGS}  
and obtain a stronger and more flexible statement with 
new parameters $\Delta$ controlling the size of character sums
and $z$ controlling the range where these sums are considered.

\begin{theorem}\label{thm:CharSmooth}
There exist absolute constants $C_0, c_0>0$ such that for  any fixed $\kappa > 0$ 
and real 
$x$, $y$, $z$, $Q$ and $\Delta$ that  satisfy
$$
(\log x)^{C_0} \le y \le z, \quad Q \le \min \{y^{c_0}, \exp(c_0 \log z/\log \log z)\}
$$
and 
$$
 \Delta \ge \max\{z^{-1}, y^{-\kappa}\}, 
$$
and an arbitrary sequence $\sA = \{a_n\}_{n=1}^\infty$ 
of complex numbers with 
$$
|a_n| \le 1, \qquad n =1, 2, \ldots, 
$$
for 
all but at most $$
E = O\(\Delta^{-2}  ( \log (1/\Delta))^2 \log x \)
$$
pairs $(q, \chi)$ with a positive integer $q \le Q$ and $\chi \in \cX_q^*$, 
for every $t \in [z,x]$ we have 
$$
S_q(\chi;\sA,t,y) = O\( \Delta \psi(t,y)\) ,
$$
where the implied constants depend only on $\kappa$. 
\end{theorem}

For example,  using that $u = O\(\log x/\log \log x\)$, we see from  Theorem~\ref{thm:CharSmooth}  that~\eqref{eq:small S}
fails for at most 
$$
E = O\((\log x)^{2\beta +8} (\log \log x)^{2}\)
$$
pairs $(q, \chi)$ under consideration, improving the bound~\eqref{eq:bound E} and holding 
in a broader range of parameters.

Finally, in Section~\ref{sec:comm} we also present an argument due to Adam Harper 
which shows that  his results~\cite{Harp} allow a more direct approach to the 
sums $S_q(\chi;x,y)$

\section{Preliminaries}

As usual, we use the expressions $F \ll G$, $G \gg F $ and $F=O(G)$ to
mean $|F|\leq cG$ for some constant $c>0$ which throughout the paper
may depend on the real parameter $\kappa > 0$.
%

First we recall   the following result of Hildebrand~\cite[Corollary~2]{Hild}:

\begin{lemma}
\label{lem:LS-SmoothInt}
For any fixed $\kappa>0$, if $1 \ge \Delta > \min\{x^{-1}, y^{-\kappa}\}$  then
$$
\psi(x+\Delta x,y) - \psi(x,y) \ll  \Delta \psi(x,y). 
$$
\end{lemma}


We use the following version of the classical large sieve
inequality, which has been given by Drappeau, Granville and Shao~\cite{DGS}:

\begin{lemma}
\label{lem:LS-SmoothChar}
There exist absolute constants $C_0, c_0>0$ such that for any real $x$, $y$ and $Q$  with 
$$
(\log x)^{C_0} \le y \le x \mand Q \le \min \{y^c_0, \exp(c_0\log x/\log \log x)\}
$$
and an arbitrary sequence $\sB = \{b_n\}_{n=1}^\infty$ 
of complex numbers, we have
$$
\sum_{q\le Q} \sum_{\chi \in \cX_q^*}
\left| S_q(\chi;\sB,x,y)  \right|^2 \ll \psi(x,y) 
\sum_{n \in \Psi(x,y)} |b_n|^2,
$$
\end{lemma}

Let $F_\xi(u)$ be the periodic function with period one for which
\begin{equation}
\label{eq:charfunc}
F_\xi(u) = \left\{  \begin{array}{ll}
1& \quad \hbox{if $0< u\le \xi$}; \\
0& \quad \mbox{if $\xi<u\le 1$}.
\end{array} \right.
\end{equation}

We now recall the following classical result of Vinogradov (see~\cite[Chapter~I,
Lemma~12]{Vin}):
 
 \begin{lemma}
\label{lem:Vinogr}
For any $\Delta$ such that
$$
0 < \Delta < \frac{1}{8} \mand \Delta\le
\frac{1}{2}\,\min\{\xi,1-\xi\},
$$
there is a real-valued function $f_{\Delta,\xi}(u) $ with the following
properties:
\begin{itemize}
\item $f_{\Delta,\xi}(u)  $ is periodic with period one;

\item $0 \le f_{\Delta,\xi}(u)  \le 1$ for all $u\in\R$;

\item $f_{\Delta,\xi}(u)  = F_\xi(u) $ if $\Delta\le u\le
\xi-\Delta$ or $\xi+\Delta\le u\le 1-\Delta$;

\item $f_{\Delta,\xi}(u) $ can be represented as a Fourier series
$$
f_{\Delta,\xi}(u)  = \xi + \sum_{j=1}^\infty \(\,g_j\,\e(ju) + h_j
\,\e(-ju)\),
$$
where $\e(u)  = \exp(2 \pi i u)$ and the coefficients $g_j,h_j$ satisfy the uniform bound
$$
\max\{|g_j|, |h_j|\}\ll\min\{j^{-1},j^{-2}\Delta^{-1}\}, \qquad j =1,2, \ldots. 
$$
\end{itemize}
\end{lemma}

\section{Proof of Theorem~\ref{thm:CharSmooth}}
\label{sec:bound char}

We cover the interval $[z,x]$ by $M = O(\log x)$ (possibly overlapping) dyadic intervals 
of the form $[m, 2m]\subseteq [z,x]$ with an integer $m$.

We now fix one of these intervals and estimate the number $E_m$ of pairs
$(q, \chi)$ with a positive integer $q \le Q$ and $\chi \in \cX_q^*$
 such that there exists $t \in [m,2m]$ for which we have 
\begin{equation}
\label{eq: large St}
\left|S_q(\chi;\sA,t,y)\right|   > c_\kappa \Delta \psi(t,y) , 
\end{equation}
where $c_\kappa> 0$ is some constant, which may depend on $\kappa$, to be chosen later. 

First,  we note that by  Lemma~\ref{lem:LS-SmoothInt}  we have
\begin{equation}
\label{eq: equal psi}
 \psi(m,y) \le  \psi(t,y)  \le  \psi(2m,y) \ll  \psi(m,y)
\end{equation}
which we use throughout the proof. 

We also note that for the number $E_{1,m}$ of such pairs
$(q, \chi)$ for which 
$$
\left|S_q(\chi;\sA,m,y)\right|  >   \Delta \psi(m,y),
$$
by Lemma~\ref{lem:LS-SmoothInt} we have 
$$
E_{1,m} \( \Delta \psi(m,y) \)^2 \ll  \psi(x,y)^2. 
$$
Hence 
\begin{equation}
\label{eq:E1}
E_{1,m} \le \Delta^{-2}.
\end{equation}
Similarly for the number $E_{2,m}$ of  pairs
$(q, \chi)$ for which 
$$
\left|S_q(\chi;\sA,2m,y)\right|   >     \Delta \psi(2m,y),
$$
we have
\begin{equation}
\label{eq:E2}
E_{2,m} \le \Delta^{-2}.
\end{equation}

So, removing these pairs $(q, \chi)$, we can now assume that
\begin{equation}
\label{eq:Sm}
\left|S_q(\chi;\sA,m,y)\right|,    \left|S_q(\chi;\sA,2m,y)\right| \le     \Delta \psi(m,y),
\end{equation}
holds. 

Let $\xi = t/m$. Using the function $F_\xi$ as in~\eqref{eq:charfunc}, we write
\begin{equation}
\begin{split}
\label{eq:Split}
S_q(\chi;\sA,t,y) &= \sum_{n \in \Psi(m,y)} a_n \chi(n)  \\
& \qquad \qquad +
\sum_{\substack{ k \ge 1\\ m+k \in \Psi(2m,y)} }a_{m+k} \chi(m+k) F_\xi(k/m).
\end{split}
\end{equation}

We can certainly assume that $\Delta < 1/8$ as otherwise the result is trivial. 

If 
\begin{equation}
\label{eq:VinogCond}
\Delta\le   \frac{1}{2}\,\min\{\xi,1-\xi\}
\end{equation}
then obviously either 
$$
\left|S_q(\chi;\sA,t,y) - S_q(\chi;\sA,m,y)\right|  \le \psi(m+ 2\Delta m,y ) - \psi(m,y ) 
$$
or 
$$
\left|S_q(\chi;\sA,2m,y) - S_q(\chi;\sA,t,y)\right|  \le \psi(2m,y ) -  \psi(2m-2\Delta m,y ) .
$$
Hence, applying Lemma~\ref{lem:LS-SmoothInt}, which is possible due to 
the  condition on $\Delta$, and  recalling~\eqref{eq:Sm}, we obtain that in this 
case 
$$
S_q(\chi;\sA,t,y) \ll \Delta  \psi(t,y ) .
$$

Therefore we can assume that~\eqref{eq:VinogCond} holds and thus
Lemma~\ref{lem:Vinogr} applies with this $\xi$ and $\Delta$. 

Considering the contribution to the above sum from $k$ with 
$$
k/m \in [0, \Delta] \cup[\xi-\Delta, \xi+\Delta] \cup [1-\Delta, 1]
$$
(that is, with $k/m$ in one of the intervals where $ F_\xi(u)$ and  $f_{\Delta,\xi}(u)$ may disagree),
and 
recalling that $|a_n|\le 1$, 
we obtain 
\begin{align*}
\bigl|\sum_{\substack{ k \ge 1\\ m+k \in \Psi(2m,y)} }&a_{m+k}  \chi(m+k) F_\xi(k/m)  \\
&-\sum_{\substack{ k \ge 1\\ m+k \in \Psi(2m,y)} }a_{m+k} \chi(m+k)f_{\Delta,\xi}(k/m)\Bigr||\\
&\qquad\qquad \qquad  \quad \   \le 
\(\psi(m+ \Delta m,y ) - \psi(m,y ) \)\\
&\qquad\qquad \qquad  \qquad \quad   + \(\psi(t+\Delta m, y)-\psi(t-\Delta m,y)\)\\
& \qquad\qquad \qquad  \qquad \qquad\quad    + 
 \(\psi(2m, y)-\psi(2m-\Delta m,y)\).
\end{align*}
Hence, applying  Lemma~\ref{lem:LS-SmoothInt} and recalling~\eqref{eq: equal psi}, we obtain 
\begin{align*}
\sum_{\substack{ k \ge 1\\ m+k \in \Psi(2m,y)} } &a_{m+k}  \chi(m+k) F_\xi(k/m)  \\
&-\sum_{\substack{ k \ge 1\\ m+k \in \Psi(2m,y)} } a_{m+k}  \chi(m+k)f_{\Delta,\xi}(k/m)\\
&\qquad\qquad  \qquad \ll \Delta   \psi(2m,y)   \ll \Delta   \psi(m,y)  \le   \Delta   \psi(t,y)  , 
\end{align*}
which  together  with~\eqref{eq:Sm} and~\eqref{eq:Split} implies
\begin{equation}
\label{eq:Sapprox1}
\begin{split}
S_q(\chi;\sA,t,y) =  
\sum_{\substack{ k \ge 1\\ m+k \in \Psi(2m,y)} } a_{m+k}  \chi(m+k)&f_{\Delta,\xi}(k/m) \\
&\quad+ O(\Delta   \psi(t,y)  ).
\end{split}
\end{equation}
Furthermore, defining $J= \rf{\Delta^{-2}}$, we see from the properties 
of the coefficients of  $f_{\Delta,\xi}(u)$   that we can approximate it  as 
$$f_{\Delta,\xi}(u) = \widetilde f_{\Delta,\xi}(u)  + O(\Delta^{-1})
$$
by a finite trigonometric polynomial
$$
\widetilde f_{\Delta,\xi}(u)  =   \sum_{j=0}^J\(\,g_j\,\e(ju) + h_j
\,\e(-ju)\), 
$$
where  for convenience  we have also defined   $g_0=h_0 = \xi/2$.
Hence, using~\eqref{eq: equal psi}, we can rewrite~\eqref{eq:Sapprox1} as
\begin{equation}
\label{eq:Sapprox2}
S_q(\chi;\sA,t,y) =  
\sum_{\substack{ k \ge 1\\ m+k \in \Psi(2m,y)} }a_{m+k}  \chi(m+k)\widetilde f_{\Delta,\xi}(k/m)   + O(\Delta   \psi(t,y)  ).
\end{equation}

Expanding the function $\widetilde f_{\Delta,\xi}(u)$ and changing the order of summation, 
we obtain 
\begin{align*}
\sum_{\substack{ k \ge 1\\ m+k \in \Psi(2m,y)} } &a_{m+k}  \chi(m+k)\widetilde f_{\Delta,\xi}(k/m) \\
&  =
 \sum_{j=0}^J\biggl(\,g_j \sum_{\substack{ k \ge 1\\ m+k \in \Psi(2m,y)} }a_{m+k}  \chi(m+k) \e(jk/m) \\
 & \qquad\qquad \qquad+ h_j
\sum_{\substack{ k \ge 1\\ m+k \in \Psi(2m,y)} }a_{m+k}  \chi(m+k) \e(-jk/m)\biggr). 
\end{align*}
Since $g_j, h_j\ll j^{-1}$ we obtain 
\begin{equation}
\label{eq:Sjm}
\sum_{\substack{ k \ge 1\\ m+k \in \Psi(2m,y)} } a_{m+k}  \chi(m+k)\widetilde f_{\Delta,\xi}(k/m)  \ll 
 \fS_{q,j}(\chi,m) , 
\end{equation}
where 
$$
 \fS_{q}(\chi,m)  =  \sum_{j=0}^J \frac{1}{j} \left|T_{q,j}(\chi,m)\right|
 $$
 with
\begin{align*}
T_{q,j}(\sA, \chi,m) & = 
\sum_{\substack{ k \ge 1\\ m+k \in \Psi(2m,y)} }a_{m+k}  \chi(m+k) \e(jk/m) \\
& = \sum_{\substack{ n > m \\ n \in \Psi(2m,y)} }a_{n}  \chi(n) \e(j(n-m)/m)   .
\end{align*}
Combining~\eqref{eq:Sapprox2} and~\eqref{eq:Sjm} together, we see that for some 
constant $c_\kappa$, depending only on $\kappa$, we have 
\begin{equation}
\label{eq:Sapprox3}
S_q(\chi;\sA,t,y) \le \frac{1}{2} c_\kappa\( \fS_{q,j}(\chi,m) +  \Delta   \psi(t,y)  \), 
\end{equation}
which is the same constant which we also use in~\eqref{eq: large St}.

We note the most  crucial for our argument point that the sums $T_{q,j}(\sA,\chi,m)$ and thus
 $\fS_{q}(\sA,\chi,m)$ 
do not depend on $t$. In particular, it is enough to estimate 
the number $E_{0,m}$  of pairs $(q, \chi)$ with a positive integer $q \le Q$ and $\chi \in \cX_q^*$  
with
\begin{equation}
\label{eq: large fSm}
 \fS_{q}(\sA,\chi,m)  \ge   \Delta \psi(m,y). 
\end{equation}

Writing $j^{-1} = j^{-1/2} \cdot j^{-1/2}$ and using the Cauchy inequality, we obtain
$$
\fS_{q}(\sA,\chi,m)^2  \ll    \sum_{j=0}^J \frac{1}{j} \left|T_{q,j}(\sA,\chi,m)\right|^2 \log J.
$$
We now apply Lemma~\ref{lem:LS-SmoothChar} for every $j =0, \ldots, J$, with the sequence 
$\sB = \{b_n\}_{n=1}^\infty$   supported only on $y$-smooth integers $n \in [m+1,2m]$
 in which case we set $b_n = a_{n}\e(j(n-m)/m)$. 
This yields the bound
\begin{align*}
\sum_{q\le Q} &\sum_{\chi \in \cX_q^*}
\fS_{q,j}(\sA,\chi,m)^2 \\
& \ll     \sum_{j=0}^J \frac{1}{j} 
\sum_{q\le Q} \sum_{\chi \in \cX_q^*}\left|T_{q,j}(\sA,\chi,m)\right|^2  \log J 
\ll \psi(2m,y)^2  ( \log J)^2, 
\end{align*}
which implies that
$$
E_{0,m} \le \Delta^{-2}  ( \log J)^2 \ll  \Delta^{-2}  ( \log (1/\Delta))^2.
$$
Therefore, recalling~\eqref{eq:E1}
and~\eqref{eq:E2} and using~\eqref{eq:Sapprox3}, we obtain 
$$E_m \le E_{0,m} + E_{1,m}+ E_{2,m}\ll \Delta^{-2}  ( \log (1/\Delta))^2
$$
for each of  $M \ll \log x$ relevant values of $m$. 
The result now follows. 

\section{Comments}
\label{sec:comm}

Here we show that the sums $S_q(\chi;x,y)$ without weights  admit a more efficient  treatement 
directly from a result of Harper~\cite[Proposition~1]{Harp}.

First we note that,  without loss of generality, in the 
hypotheses  of Theorem~\ref{thm:CharSmooth}, we can assume that $\kappa$ is sufficiently small
(since otherwise, adjusting the value of the constant $c_0$, 
we can use  the trivial bound $Q^2$ for the number of ``bad''  pairs $(q,\chi)$).

We now set 
$$
D = \Delta^{-1} \log x, \qquad H = D^{50}, \qquad \varepsilon =\frac{C}{\log y} + \frac{4\log D}{\log z},
$$ 
where $C$ is the absolute constant  of~\cite[Proposition~1]{Harp}.

It is easy to check  that all the assumptions of~\cite[Proposition~1]{Harp}, 
used with $t$ in place of $x$, are satisfied for any such character and any  
$t$ with $z \leq t \leq x$. Indeed, we first note that by~\cite[Lemma~2]{HiTe} 
the parameter $\alpha(t,y)$ satisfies $\alpha(t,y) = 1 + o(1)$.  Hence, for a
sufficiently large $C_0$ and a sufficiently small $\kappa$ we have 
\begin{equation}
\begin{split}
\label{eq: cond1}
\frac{C}{\log y}  < \varepsilon &\le \frac{C}{\log y} + \frac{4 \log (1/\Delta) + 4 \log \log x}{\log z} \\
&\le \frac{C}{\log y} + \frac{4\kappa \log y+ 4 \log \log x}{\log z} \\
&\le \frac{C}{\log y} +  4\kappa + \frac{4}{C_0} < \alpha(t,y)/2,
\end{split}
\end{equation}
(provided that $x$ is large enough). 
We also note that 
$$
y^ \varepsilon = C  \exp\( \frac{4\log D \log y}{\log z}\) \ll D^4.
$$

We now assume that $C_0$ is sufficiently large and $\kappa$ is sufficiently small, 
so that  we have
$$
\kappa + 1/C_0< d/50,
$$
where $d$ is the absolute constant  of~\cite[Proposition~1]{Harp}. 
Hence, recalling the inequalities  $(\log x)^{C_0} \le y \le z$, 
we also have 
\begin{equation}
\begin{split}
\label{eq: cond2}
y^ {0.9\varepsilon} (\log t)^2 \le   D^{3.6}  &(\log x)^2 \le D^{5.6}  \le H \\
& \quad \le \Delta^{-50}   (\log x)^{50} 
\le y^{50 \kappa+50/C_0} \le    z^d \le t^d.
\end{split}
\end{equation}
Finally, we also verify that  for $q \le Q$ and an appropriate choice of 
the constants $c_0$,  $C_0$ and $\kappa$ we have 
\begin{equation}
\label{eq: cond3}
(Hq)^A \le \(Q (\Delta^{-1} \log x)^{50}\)^A \le (y^{c_0 +  50 \kappa + 50/C_0})^A  \le y.
\end{equation}
Thus the inequalities~\eqref{eq: cond1}, \eqref{eq: cond2} and~\eqref{eq: cond3}
ensure that~\cite[Proposition~1]{Harp} applies and implies  that 
for any modulus $q \le Q$ and a primitive character $\chi$ modulo $q$, such that the $L$-function 
$L(s,\chi)$ has no zeros in the region 
$$
\Re(s) > 1-\varepsilon \mand |\Im(s)| \leq H,
$$ 
we have
\begin{align*}
 \Psi(t,y;\chi) &\ll \Psi(t,y) \sqrt{\log t \log y} (t^{-0.3\varepsilon}\log H + H^{-0.02})\\
                     &\leq \Psi(t,y) \sqrt{\log t \log y} (D^{-1.2}\log H + D^{-1})\\
                      &  \ll \Delta \Psi(t,y)
\end{align*}
for all $z \leq t \leq x$. 

It only remains to see how many pairs $(q,\chi)$ do not have such a zero-free region. 
By the zero-density estimates of Huxley~\cite{Hux} and Jutila~\cite{Jut}
(see also~\cite[Section~2]{Harp}, where such results are conveniently 
summarised),  this number is at most $(Q^2H)^{(12/5 + o(1))\varepsilon}$. 
First we note that $D  \le y^{\kappa + 1/C_0}$, which implies
$$
H^\varepsilon = \exp\( \frac{200(\log D)^2}{\log z} + O\(\frac{\log D}{\log y}\) \)
 \ll   \exp\( \frac{200(\log D)^2}{\log z} \) .
$$
It is also useful to note that since  $Q \le y^{c_0}$,  we have
$$
Q^\varepsilon  = Q^{4\log D/\log z +  O(1/\log y)}  \ll Q^{4\log D/\log z}.
$$

Recalling that $Q$ is at most a small power of $y$ we now derive the number $E$ of ``bad''  pairs $(q,\chi)$
is at most 
$$
E \ll Q^{(96/5 + o(1))\log D/\log z} \exp\((480+o(1)) \frac{\(\log D\)^2}{\log z}\).
$$
with the above choice of $D$. 
In particular, if 
\begin{equation}
\label{eq: Delta z}
\Delta \ge \exp\(-c_1 \sqrt{\log z}\) \mand z \ge \exp\(c_2(\log \log x)^2\)
\end{equation}
for some absolute constants $c_1,c_2 >0$, 
then, using the condition $Q  \le \exp(c_0 \log z/\log \log z)$ we obtain 
$$
 Q^{\log D/\log z} \le \exp(c_0 \log D/\log \log z) = D^{O(1/ \log \log z)}
$$
and 
$$
\exp\(\frac{\(\log D\)^2}{\log z}\) = \exp\( O\( \frac{\(\log(1/\Delta) + \log \log x\)^2 }{\log z}\) \)= O(1).
$$
Hence
$$
E \ll D^{O(1/ \log \log z)}.
$$
Thus, if in instead of~\eqref{eq: Delta z} we impose more stringent conditions
$$
\Delta \ge \frac{c_1}{\log z} \mand   z \ge  \exp\( (\log  x)^{c_2}\)
$$ 
then we see that $E = O(1)$.

Finally, we  note that~\cite[Theorem~3]{Harp} can also be used to estimate the sums
$S_q(\chi;x,y)$. Furthermore, the sums $S_q(\chi;\sA, x,y)$ with weights given by multiplicative 
functions, such as the {\it M{\"o}bius function\/}, can be treated  within the same technique.
However,  this approach  does not seem to apply to the sums
$S_q(\chi;\sA, x,y)$ with arbitrary  weights. 

\section*{Acknowledgement}

The author is grateful to Adam Harper for proving the argument used in
Section~\ref{sec:comm}  together with the generous permission to present it here
as well as further comments. 

Some parts of this work were done when
the author was visiting 
the Max Planck Institute for Mathematics, Bonn,
and the Institut de Math{\'e}matiques de Jussieu, Universit{\'e} Paris Diderot, 
whose  generous 
support and hospitality are gratefully acknowledged. 

This work was also partially supported    
by the Australian Research Council Grant~DP140100118.

\end{document}